\documentclass[review,3p]{elsarticle}

\usepackage{graphics} 
\usepackage{epsfig} 
\usepackage{mathptmx} 
\usepackage{times} 
\usepackage{amsmath} 
\usepackage{amssymb}  
\usepackage{psfrag}
\usepackage{textcomp}
\usepackage{amsthm}
\usepackage{setspace}
\usepackage{wasysym}
\usepackage{latexsym}
\usepackage{xfrac}
\usepackage{relsize}
\usepackage{pmat}
\usepackage{url}
\usepackage{comment}
\usepackage{pdfsync}

\usepackage[cal=cm]{mathalfa} 

\newcommand{\aHalf}{\frac{1}{2}}

\DeclareMathOperator*{\argmin}{arg\;min}
\DeclareMathOperator{\co}{co}
\newcommand{\del}{\nabla}

\newcommand{\PXi}{\mathbf{P}_{\Xi}}

\newcommand{\Pxst}{\mathbf{P}_{X^*}}

\newcommand{\Pa}{\mathbf{P}_{\!\!\mathcal{A}}}

\newcommand{\Bgo}[1]{\bar{B}_{#1}(\Gamma_0)}

 \theoremstyle{remark} 
\newtheorem{assump}{A\!\!}[section]

\newtheorem{lemma}{Lemma}[section]
\newtheorem{corollary}{Corollary}[section]

\newtheorem{example}{Example}[section]
\newtheorem{theorem}{Theorem}[section]
\newtheorem{defn}{Definition}[section]
\newtheorem{rem}{Remark}[section]

\usepackage{hyperref}

\begin{document}

\begin{frontmatter}

\title{A Semiglobal, Practical, Strict Pseudogradient Property for Iterative Methods}

\author[myaddress]{Karla Kvaternik\corref{CorrAuth}}

\cortext[CorrAuth]{This research was supported by the Army Research Office Grant W911NF-14-1-0431. Correspondence can be sent to: \ead{karlak@princeton.edu}}

\address[myaddress]{Mechanical and Aerospace Engineering, Engineering Quadrangle, Princeton University, Princeton, NJ.}

\begin{abstract}
We consider a class of iterative numerical methods and introduce the notion of \emph{semiglobally, practically, strictly pseudogradient} (SPSP) search directions. We demonstrate the relevance of the SPSP property in modelling a variety of optimization algorithms, including those subject to absolute and relative errors. We show that the attractors of iterative methods with SPSP search directions exhibit semiglobal, practical, asymptotic stability. Moreover, the SPSP property is \emph{robust} in the sense that perturbations of SPSP search directions also have the SPSP property. 
\end{abstract}

\begin{keyword}
discrete-time nonlinear systems\sep iterative numerical methods\sep semiglobal, practical, asymptotic stability\sep robustness\sep Lyapunov analysis
\end{keyword}

\end{frontmatter}


%
%


\section{Introduction}
We consider a class of iterative numerical methods of the form
\begin{equation}\label{eq:iterative}
y = \PXi\big[y - \alpha s\big], \quad y \in\mathbb{R}^n, \; s\in\Psi(y),
\end{equation}
where $\alpha\in\mathbb{R}_{++}$ is a tunable gain, $\Psi:\mathbb{R}^n\rightrightarrows\mathbb{R}^n$ is a set of \emph{search directions} at each $y$, and $\mathbf{P}_{\Xi}(\cdot)$ is the orthogonal projection onto a closed, convex set $\Xi\subset\mathbb{R}^n$.

Although \eqref{eq:iterative} need not be associated with optimization of some objective, it models a wide variety of convex optimization methods designed to iteratively approximate the set of minimizers
\begin{equation}\label{eq:optProb}
X^* = \argmin_{y\in\Xi} J(y)
\end{equation}
of a function $J:\mathbb{R}^n\rightarrow\mathbb{R}$, over a feasible set $\Xi$.

We make three contributions in this manuscript. First, we introduce the notion of \emph{semiglobally, practically, strictly pseudogradient} (SPSP) search directions, characterized by an inequality involving elements of $\Psi(\cdot)$ and the gradient of a surrogate (i.e. Lyapunov) function $V(\cdot)$. This notion extends those of \emph{pseudogradient}, and \emph{strongly pseudogradient} search directions found in \cite{Polyak87}, for example. 

A primary motivation for introducing the SPSP property is that the qualitative behavior of sequences generated by numerical methods of the form \eqref{eq:iterative} with SPSP search directions is guaranteed to exhibit not only attractivity, but also stability with respect to a posited attractor $\mathcal{A}\subset\mathbb{R}^n$ for \eqref{eq:iterative}. Specifically, in $\S$\ref{sec:SPSP=SPAS} we show that the class of systems \eqref{eq:iterative} having SPSP search directions have a attractor which is \emph{semiglobally, practically, asymptotically stable} (SPAS), and this demonstration constitutes our second contribution in this manuscript. 

We emphasize SPAS because the analytical focus in the optimization literature tends to center almost exclusively on the concept of attractivity and the quantification of convergence rates, while stability is almost never explicitly addressed. Meanwhile, the modern proliferation of data-driven, ``cyber-physical systems" control applications may involve the dynamic interaction of optimization algorithms like \eqref{eq:iterative} and real physical systems with dynamics of their own \cite{KvaternikPHD}. In such applications, the concept of stability often relates to safety and reliability, and is therefore fundamentally important. 

Another motivation for introducing the SPSP property is its generality; even within the context of optimization, there is a wide variety of algorithms that satisfy this property, including subgradient methods, quasi-Newton methods, weighted gradient methods, and methods affected by absolute and relative deterministic errors (q.v. $\S$\ref{sec:Examples}). In working with the generic form \eqref{eq:iterative} and the SPSP property, we are able to draw conclusions that apply broadly to all such algorithms. 

A third motivation for appealing to the SPSP property in the analysis of iterative numerical algorithms is that methods  \eqref{eq:iterative} having SPSP search directions are \emph{robust}; in $\S$\ref{sec:robustness}, we show that if $\Psi(\cdot)$ is SPSP with respect to a function $V(\cdot)$, then a perturbation of $\Psi(\cdot)$ involving absolute and relative errors is also SPSP with respect to the same function $V(\cdot)$, provided the errors are sufficiently small. This demonstration constitutes our third contribution in this manuscript.

Finally, we make the observation that in the context of optimization, convexity of neither the objective function $J(\cdot)$ nor the surrogate function $V(\cdot)$ is required for SPSP to hold. Thus, the SPSP property may potentially help guide the development of more broadly applicable optimization methods.  

A minor contribution of this manuscript is a set of analytic tools (Lemmas \ref{lem:sublevelSetContainment} and \ref{lem:underEstimation}) for establishing the SPSP property in a broad variety of contexts.

The work in this manuscript draws  conceptually on a combination of the Lyapunov-based techniques used to study optimization algorithms in \cite{Polyak87}, and the notion of SPAS, which appears to have been initially formalized in \cite{TeelSCL99}, \cite{TeelCSM} (c.f. \cite{KvaternikSPAS}). Also related are the results found in \cite{TeelCDC00a}. The analytic tools and robustness characterizations that we provide here can be regarded as alternative to those in \cite{TeelCDC00a}.

\paragraph{Notation and Preliminaries}\label{sec:prelim}
All vector norms $\|\cdot\|$ are Euclidean. If $S\subset\mathbb{R}^n$ is closed and $x_o\in\mathbb{R}^n$, then
\begin{equation}\label{eq:proj}
\mathbf{P}_S(x_o) = \argmin_{x\in S} \|x_o-x\|
\end{equation}
is the orthogonal projection of $x_o$ onto $S$. The set of non-negative real numbers is denoted by $\mathbb{R}_+$, while the set of positive real numbers is denoted by $\mathbb{R}_{++}$. For a point $x_o\in\mathbb{R}^n$, and $r\in\mathbb{R}_{++}$, $\bar{B}_r(x_o)=\{x\in\mathbb{R}^n: \|x-x_o\|\leq r\}$ and $B_r(x_o)=\{x\in\mathbb{R}^n: \|x-x_o\|<r\}$. For a compact set $S\subset\mathbb{R}^n$, $\bar{B}_{r}(S)=\{x\in\mathbb{R}^n\;\big|\; \|x-\mathbf{P}_S(x)\|\leq r\}$, while $B_{r}(S)=\{x\in\mathbb{R}^n\;\big|\; \|x-\mathbf{P}_S(x)\|< r\}$. Strict containment of a set $A$ inside a set $B$ is denoted as $B\supset A$, and means that for any  $y\in A$ it is possible to find an $\epsilon>0$ such that $B_{\epsilon}(y)\subset B$. Set containment is denoted as $B\supseteq A$. The set of continuous (resp. continuously differentiable) functions from $\mathbb{R}^n$ into $\mathbb{R}^m$ is denoted by $C^0[\mathbb{R}^n,\mathbb{R}^m]$ (resp. $C^1[\mathbb{R}^n,\mathbb{R}^m]$). We say that a function $V\in C[\mathbb{R}^n,\mathbb{R}_+]$ is \emph{positive definite} with respect to a closed set $S$ on a set $\Omega\supset S$ if, $V(S)=\{0\}$, and $V(x)>0$ for all $x\in\Omega\backslash S$. A function $V:x\mapsto V(x)$ on $\mathbb{R}^n$ is \emph{radially unbounded} with respect to a closed set $S\subset\mathbb{R}^n$, if for any $B\in\mathbb{R}$, there exists an $r\in\mathbb{R}_{++}$ such that $V(x)>B$, for all $x\in\mathbb{R}^n\backslash \bar{B}_r(S)$. The gradient of a differentiable function $J:\mathbb{R}^n\rightarrow\mathbb{R}$ is denoted by $\del J(\cdot)$. The subdifferential of a convex function $J:\mathbb{R}^n\rightarrow\mathbb{R}$ at a point $x\in\mathbb{R}^n$ is denoted by $\partial J(x)$, and consists of all vectors $g\in\mathbb{R}^n$ such that ${J(x+y)-J(x)\geq f(x) + g^Ty}$ holds for all $y\in\mathbb{R}^n$. Using the definition of Frechet differentiability, it can be shown that for $V(y) = \tfrac{1}{2}\|y-\Pxst(y)\|^2$, $\del V(y) = (y-\Pxst (y))$. For a sequence $(x(t))_{t=0}^{\infty}$ we use $x$ to stand for $x(t)$, and $x^+$ to stand for $x(t+1)$. For $S\subset\mathbb{R}^n$, $\co(S)$ is its convex hull.

\section{Semiglobally, Practically, Strictly Pseudogradient Search Directions}\label{sec:SPSP}

In the following definitions, we specify a property of $\Psi(\cdot)$ that allows us to draw general conclusions about the qualitative behavior of \eqref{eq:iterative}. This property is a generalization of the  \emph{pseudogradient} property discussed in \cite{Polyak87} (q.v. Remark \ref{rem:pseudograd1}) and because of its link to semiglobal, practical, asymptotic stability (q.v. $\S$\ref{sec:SPSP=SPAS}), we refer to it as the \emph{semiglobal, practical, strict pseudogradient} (SPSP) property. We motivate our definition of SPSP search directions with several examples in $\S$\ref{sec:Examples}.

\begin{defn}[Semiglobally,  practically strictly pseudogradient search directions]\label{def:SPStrictPseudo}
Consider a multi-valued vector field
$\Psi:\mathbb{R}^n\rightrightarrows \mathbb{R}^n$ and a differentiable function $V:\mathbb{R}^n\rightarrow\mathbb{R}_+$, which is positive definite and radially unbounded with respect to a compact set $\mathcal{A}\subset\mathbb{R}^n$. We say that $\Psi(\cdot)$ is  \emph{semiglobally, practically, strictly pseudogradient} (SPSP) with respect to $V(\cdot)$ on a set $\Xi\supset \mathcal{A}$, if for some $\epsilon\in\mathbb{R}_+$ and $b\in\mathbb{R}_+$, and for any $\sigma\in\mathbb{R}_{++}$  (with $\sigma >\epsilon$), 
\begin{equation*}
\del V(y)^Ts\geq -b,\quad \forall y\in\Xi\cap\bar{B}_{\epsilon}(\mathcal{A}),\;\forall s\in\Psi(y),
\end{equation*}
and there exists a function $\phi_{\sigma,\epsilon}\in C^0[\mathbb{R}^n,\mathbb{R}]$ which is positive on $\Xi\cap \big(\bar{B}_{\sigma}(\mathcal{A})\backslash B_{\epsilon}(\mathcal{A})\big)$ and radially unbounded with respect to $\mathcal{A}$ on $\Xi$, such that 
\begin{equation*}
\del V(y)^Ts\geq \phi_{\sigma,\epsilon}(y), \quad \forall y\in\Xi\cap (\bar{B}_{\sigma}(\mathcal{A})\backslash B_{\epsilon}(\mathcal{A})),\;\forall s\in\Psi(y).
\end{equation*}

If these conditions hold with $\epsilon=b=0$, we say that $\Psi(\cdot)$ is \emph{semiglobally, strictly pseudogradient} (SSP) with respect to $V(\cdot)$ on $\Xi$. 

If these conditions hold with $\bar{B}_{\sigma}(\mathcal{A})=\mathbb{R}^n$, and independently of $\sigma$, then we say that $\Psi(\cdot)$ is \emph{(globally) practically, strictly pseudogradient} (PSP) with respect to $V(\cdot)$ on $\Xi$. 

If these conditions hold with $\bar{B}_{\sigma}(\mathcal{A})=\mathbb{R}^n$ and $\epsilon = b = 0$, we say that $\Psi(\cdot)$ is \emph{strictly pseudogradient} (SP) with respect to $V(\cdot)$ on $\Xi$. 
$\diamondsuit$\end{defn}

\begin{rem}
Definition \ref{def:SPStrictPseudo} is a direct generalization of Definition  3.3.1 in \cite{KvaternikPHD}.
$\diamondsuit$\end{rem}

\begin{rem}
When algorithm \eqref{eq:iterative} pertains to the optimization problem \eqref{eq:optProb}, it is typically the case that $\mathcal{A}\equiv X^*$. 
$\diamondsuit$\end{rem}

\begin{rem}[The Rationale for Definition \ref{def:SPStrictPseudo}, and its relationship to similar notions]\label{rem:pseudograd1}
Definition \ref{def:SPStrictPseudo} complements and generalizes the concepts of \emph{pseudogradient}, and \emph{strongly pseudogradient} search directions found in $\S$2.2, \cite{Polyak87}. 

Although the strict pseudogradient property in Definition \ref{def:SPStrictPseudo} is weaker than the \emph{strong pseudogradient property} in \cite{Polyak87}, which requires that
\begin{equation}
\del V(\xi)^Ts \geq \tau V(\xi),\quad \forall s\in\Psi(\xi),\;\forall \xi\in\mathbb{R}^n,
\end{equation}
and some $\tau>0$, it is stronger than the \emph{pseudogradient property} in \cite{Polyak87}, which requires that
\begin{equation}
\del V(\xi)^Ts \geq 0,\quad \forall s\in\Psi(\xi),\;\forall \xi\in\mathbb{R}^n.
\end{equation}
$\diamondsuit$\end{rem}

\section{Optimization Examples Motivating the Use of Definition \ref{def:SPStrictPseudo}}\label{sec:Examples}
One of the motivations for introducing the SPSP property is its generality; even within the context of optimization, there is a wide variety of algorithms that satisfy this property, and in working with the generic form \eqref{eq:iterative}, we are able to draw conclusions that apply broadly to all such algorithms. 

In this section we consider the optimization problem \eqref{eq:optProb} and we work through a number of examples involving various assumptions on $J(\cdot)$ and various subgradient-related choices of $\Psi(\cdot)$, to show how the SPSP property applies. We begin with a number of basic examples of algorithms employing SP search directions in $\S$\ref{sec:basicEx}. In $\S$\ref{sec:errorEx} we examine algorithms employing SSP, PSP and SPSP search directions involving relative and absolute deterministic errors, and in $\S$\ref{sec:nonconvexEx} we observe that gradients of non-convex functions can also satisfy Definition \ref{def:SPStrictPseudo}.

Although all our examples involve either $V(y) = \tfrac{1}{2}\|y-\Pxst(y)\|^2$ or $V(y) = J(y) - J^*$ (where $J^*:=\min_{y\in\Xi} J(y)$), these need not be considered as the only possible choices for $V(\cdot)$ in the analysis of algorithms like \eqref{eq:iterative}.

\subsection{Basic Examples}\label{sec:basicEx}
\begin{example} 
Suppose $J(\cdot)$ is convex and differentiable. Then $\Psi(y) = \{\del J(y)\}$ is trivially SP with respect to either $V(y) = J(y) - J^*$ or $V(y) = \tfrac{1}{2}\|y-\Pa(y)\|^2$ on $\mathbb{R}^n\supset \mathcal{A}\equiv X^*$. For $V(y) = J(y) - J^*$, we may take $\phi_{\sigma,\epsilon}(y) = \|\del J(y)\|^2$, while for $V(y) = \tfrac{1}{2}\|y-\Pa(y)\|^2$, we could take $\phi_{\sigma,\epsilon}(y)=J(y)-J^*$. In both cases, $\bar{B}_{\sigma}(X^*)=\mathbb{R}^n$, and $\epsilon = b=0$.
$\diamondsuit$ \end{example}

\begin{example}
Suppose $J(\cdot)$ is a convex function that is not differentiable everywhere. By Proposition 4.2.1 in \cite{Bertsekas03}, the convexity of $J(\cdot)$ suffices to guarantee that its subdifferential  $\partial J(y)$ is nonempty $\forall y\in\mathbb{R}^n$. From its definition, a subgradient $s(y)\in\partial J(y)$ satisfies
\begin{equation*}
(y-y_o)^Ts(y) \geq J(y)-J(y_o), 
\end{equation*}
for all $y_o\in\mathbb{R}^n$, and in particular for $y_o=\Pxst (y)$. 

From this, and the fact that $\del \tfrac{1}{2}\|y-\Pa(y)\|^2=(y-\Pa(y))$, we observe that $\Psi(y)=\partial J(\cdot)$ is SP with respect to $V(y)=\tfrac{1}{2}\|y-\Pa(y)\|^2$ on $\mathbb{R}^n\supset X^*$, with $\phi_{\sigma,\epsilon}(y) = J(y)-J^*$, $\bar{B}_{\sigma}(X^*)=\mathbb{R}^n$, and $\epsilon = b=0$.
$\diamondsuit$\end{example}

\begin{example}\label{ex:stronglyConvexNoErrors}
Suppose that $J(\cdot)$ is a $c$-strongly convex function that is not differentiable everywhere. In that case for $V(y) = \tfrac{1}{2}\|y-\Pxst(y)\|^2$ and $s(y)\in\partial J(y)$, we have that
\begin{align*}
\del V(y)^Ts(y) &\geq J(y)-J(y_o)
\\
& \geq \tfrac{c}{2}\|y-\Pxst(y)\|^2
\\
& = \tfrac{c}{2} V(y),
\end{align*}
and we conclude that $\partial J(y)$ is SP with respect to $V(y) = \tfrac{1}{2}\|y-\Pxst(y)\|^2$ on $\mathbb{R}^n\supset X^*$, with $\phi_{\sigma,\epsilon}(y) = V(y)$, $\bar{B}_{\sigma}(X^*)=\mathbb{R}^n$, and $\epsilon = b=0$. This case corresponds to the notion of ``strongly pseudogradient'' search directions discussed in $\S$2.2.3 of \cite{Polyak87}. 
$\diamondsuit$ \end{example}

\begin{example}[Weighted Gradient Methods]
Consider the case in which $J(\cdot)$ is differentiable, and for all $t\in\mathbb{N}$, $s\in\Psi(y(t)) = \{H(t)\del J(y(t))\}$ and $H(t)=H(t)^T>0$. This is a large class of algorithms including Newton's method, diagonally scaled gradient descent, the Gauss-Newton method for the sum of squares of functions in $\mathbb{R}$, and quasi-Newton methods (q.v. $\S$1.2 and $\S$1.7 in \cite{Bertsekas95}), for example. 

Such search directions are trivially SP with respect to $V(y) = J(y)-J^*$, when $J(\cdot)$ is strictly convex. In that case $\bar{B}_{\sigma}(X^*)=\mathbb{R}^n$, $\epsilon = b=0$, and we may take $\phi_{\sigma,\epsilon}(y) = \lambda_m\|\del J(y)\|^2$, where $\lambda_m = \min \cup_{t\in\mathbb{N}} \sigma(H(t))>0$.
$\diamondsuit$\end{example}

\subsection{Optimization Algorithms Involving Deterministic Errors}\label{sec:errorEx}
In the following examples, we consider problem \eqref{eq:optProb}, and search directions of the form
\begin{equation}\label{eq:sWithError}
s(y) = g(y) + \eta(y),
\end{equation}
where $g(\cdot)$ may represent a subgradient, or a weighted gradient of $J(\cdot)$, and $\eta(\cdot)$ represents an error with both relative and absolute deterministic components (c.f. $\S$4.1.2 in \cite{Polyak87}). Specifically, we assume that for any $\sigma\in\mathbb{R}_{++}$, there exist positive, real numbers $a$ and $r$ such that
\begin{equation}\label{eq:errModel1}
\|\eta(y)\|\leq a + r\|y-\Pxst(y)\|,
\end{equation}
for all $y\in\bar{B}_{\sigma}(X^*)$.

The following two examples demonstrate the relevance of the error model \eqref{eq:sWithError}-\eqref{eq:errModel1}.

\begin{example}[Forward Euler Approximations of $\del J(\cdot)$ incur Absolute Errors]
Consider a function $J\in C^1[\mathbb{R}^n,\mathbb{R}]$ having a locally Lipschitz gradient, and a forward-Euler approximation of its gradient along the $i$th standard basis vector $e_i$:
\begin{equation}
s(y)^Te_i = 
\frac{J(y+\mu e_i)-J(y)}{\mu},
\end{equation}
where $\mu\in\mathbb{R}_{++}$ is a tunable approximation step-size. 

By assumption, there exists a Lipschitz constant $L_{\del J}$ for $\del J(\cdot)$ associated to any compact $\Omega\subset\mathbb{R}$. Then, by rearranging the standard expression of the Mean Value Theorem (q.v. Proposition 1.1.12 in \cite{Bertsekas03}, for example), it can be shown that for any $y,\nu\in\mathbb{R}^n$ such that $y,y+\nu\in\Omega$, 
\begin{equation}\label{eq:FDbound}
\| J(y+\nu) - J(y) - \del J(y)^T\nu\| \leq L_{\del J}\|\nu\|^2.
\end{equation}
Taking $\nu = \mu e_i$, we obtain that 
\begin{equation}
s(y)^Te_i = \del J(y)^Te_i + \eta_i(y), 
\end{equation}
where 
\begin{equation}
|\eta_i(y)| \leq \mu L_{\del J}.
\end{equation}
We therefore conclude that search directions obtained as forward-difference approximations of gradients can be modelled by \eqref{eq:sWithError}-\eqref{eq:errModel1}, with $g(y)=\del J(y)$, $r=0$ and $a = \sqrt{n}\mu L_{\del J}$.
$\diamondsuit$\end{example}

\begin{example}[Relative Errors in Weighted Gradient Methods]
When $J(\cdot)$ is continuously differentiable with a locally Lipschitz gradient, weighted gradient and quasi-Newton methods take the form \eqref{eq:iterative}, with $\Psi(y) = \{H\del J(y)\}$. Computing the matrix $H$ typically introduces numerical errors so that the implemented search direction takes the form
\begin{align}
s(y) &= (H+\tilde{H})\del J(y)
\nonumber\\
& := g(y) + \eta(y).
\end{align}
By assumption, there exists a number $L_{\del J}$ such that for all $y\in\bar{B}_{\sigma}(X^*)$,
\begin{align*}
\|\eta(y)\| &= \|\tilde{H} \del J(y)\|
\\
& = \|\tilde{H} \big(\del J(y)-\del J(\Pxst(y)) \big)\|
\\
&\leq (L_{\del J} \|\tilde{H}\| ) \|y-\Pxst(y)\|,
\end{align*}
and we conclude that quantization or inversion errors in weighted gradient methods can be modelled by \eqref{eq:sWithError}-\eqref{eq:errModel1}, with $g(y)= H\del J(y)$, $r=L_{\del J} \|\tilde{H}\| $ and $a = 0$.
$\diamondsuit$\end{example}

In the following examples, we will show how to derive conditions under which a set $\Psi(\cdot)$ of search directions \eqref{eq:sWithError} is SPSP with respect to  $V(y) = \tfrac{1}{2}\|y-\Pxst(y)\|^2$, under various assumptions on $J(\cdot)$.

Assuming that $J(\cdot)$ is convex and $g(y)\in\partial J(y)$, for this $V(\cdot)$, the error model \eqref{eq:sWithError}-\eqref{eq:errModel1} yields
\begin{align}
\del V(y)^Ts(y) &= (y-\Pxst(y))^T (g(y)+\eta(y))
\nonumber\\
&\geq (J(y)-J^*) - a \|y-\Pxst(y)\| - r \|y-\Pxst(y)\|^2.
\label{eq:start}
\end{align}

\begin{example}[Linear objective]
Suppose $J(y) = c\|y - \Pxst(y)\|$, $c>0$. Then, from \eqref{eq:start} we observe that  $\Psi(\cdot)$ is SSP with respect to $V(y) = \tfrac{1}{2}\|y-\Pxst(y)\|^2$ provided $a$ and $r$ are such that $a<c$, and for any given $\sigma>0$, $r<\tfrac{c-a}{\sigma}$. In this case $\phi_{\sigma,\epsilon}(y) = r\|y-\Pxst(y)\|\big(\tfrac{c-a}{r}-\|y-\Pxst(y)\|\big)$, and $\epsilon=b=0$.
$\diamondsuit$ \end{example}

\begin{example}[Strongly convex objective]\label{ex:stronglyConvexObj}
Suppose that $J(\cdot)$ is $c$-strongly convex on $\mathbb{R}^n$ -- that is, there exists a number $c>0$ such that for all $y\in\mathbb{R}^n$, $J(y)-J^*$ is globally underestimated by the quadratic function $ \tfrac{c}{2}\|y-\Pxst(y)\|^2$. Then from \eqref{eq:start}, we have that
\begin{align}\label{eq:bound}
\del V(y)^Ts(y) &\geq \tfrac{c-2r}{2} \|y-\Pxst(y)\| \big( \|y-\Pxst(y)\| - \tfrac{2a}{c-2r}\big),
\end{align}
from which we see that $\Psi(\cdot)$ is PSP with respect to $V(y) = \tfrac{1}{2}\|y-\Pxst(y)\|^2$ provided $r<\tfrac{c}{2}$. In that case Definition \ref{def:SPStrictPseudo} is satisfied with $\phi_{\sigma,\epsilon}(\cdot)=\tfrac{c-2r}{2} \|y-\Pxst(y)\| \big( \|y-\Pxst(y)\| - \tfrac{2a}{c-2r}\big)$, $\bar{B}_{\sigma}(X^*)=\mathbb{R}^n$, $\epsilon = \tfrac{2a}{c-2r}$, and $b=\tfrac{a^2}{2(c-2r)}$. 
$\diamondsuit$ \end{example}

\begin{example}[Convex objective]\label{ex:convexObj}
Suppose $J(y)$ is convex, though not necessarily differentiable. Then, $J(y)-J^*$ is continuous, radially unbounded, and positive definite with respect to $X^*$. We may therefore apply Corollary \ref{cor:underEstimation} to conclude that for any positive, real numbers $\hat{\sigma}$ and $\hat{\epsilon}$, with $\hat{\sigma}>\hat{\epsilon}$, there exists a number $c\in\mathbb{R}_{++}$ such that 
\begin{equation}\label{eq:underestimate}
J(y)-J^*\geq \frac{c}{\hat{\sigma}^2} \|y-\Pxst(y)\|^2
\end{equation}
on the ``band'' $\bar{B}_{\hat{\sigma}}(X^*)\backslash B_{\hat{\epsilon}}(X^*)$. For our purposes here, $c$ may be constructed as
\begin{equation}\label{eq:c}
c = \min_{y\in\partial \bar{B}_{\hat{\epsilon}}(X^*)}  \; (J(y)-J^*), 
\end{equation}
so that the $c$-sublvel set of $J(y)-J^*$ is contained inside $\bar{B}_{\hat{\epsilon}}(X^*)$ (c.f. the proof of Lemma \ref{lem:sublevelSetContainment}).

Combining \eqref{eq:start} and \eqref{eq:underestimate}, we see that on $\bar{B}_{\hat{\sigma}}(X^*)\backslash B_{\hat{\epsilon}}(X^*)$,
\begin{align}\label{eq:bound2}
\del V(y)^Ts(y) &\geq \tfrac{c-\hat{\sigma}^2r}{\hat{\sigma}^2}\|y-\Pxst(y)\|^2 \big(\|y-\Pxst(y))\| - \tfrac{\hat{\sigma}^2a}{c-\hat{\sigma}^2r}\big).
\end{align}

We claim that $\Psi(\cdot)$ is SPSP with respect to $V(y) = \tfrac{1}{2}\|y-\Pxst(y)\|^2$, provided that for any given $\sigma$, the error model \eqref{eq:sWithError}-\eqref{eq:errModel1} satisfies the following conditions:
\begin{align}
r &<\tfrac{c}{\sigma^2},
 \label{eq:rbound}
\\
a &<\tfrac{c-\sigma^2r}{\sigma}. 
\label{eq:abound}
\end{align}
To see this,  in \eqref{eq:bound2} we assign $\hat{\sigma}=\sigma$ and $\hat{\epsilon}=\epsilon$, where 
\begin{equation}
\epsilon:=\tfrac{\sigma^2 a}{c-\sigma^2r}.
\end{equation}
For this choice of $\epsilon$, \eqref{eq:abound} ensures that $\sigma>\epsilon$, as required. On the other hand, \eqref{eq:rbound} ensures that 
\begin{equation}
\phi_{\sigma,\epsilon}(y) := \tfrac{c-\sigma^2r}{\sigma^2}\|y-\Pxst(y)\|^2 \big(\|y-\Pxst(y)\| - \tfrac{\sigma^2a}{c-\sigma^2r}\big)
\end{equation}
is positive on $\bar{B}_{\sigma}(X^*)\backslash B_{\epsilon}(X^*)$. Finally, from \eqref{eq:start} we observe that for all $y\in\bar{B}_{\epsilon}(X^*)$, $\del V(y)^Ts(y)\geq -b$, where
\begin{align*}
b &= a\epsilon + r\epsilon^2
\\
&= a\big(\tfrac{\sigma^2 a}{c-\sigma^2r}\big) + r\big(\tfrac{\sigma^2 a}{c-\sigma^2r}\big)^2,
\end{align*}
completing the proof of our claim. 
$\diamondsuit$ \end{example}


\subsection{Non-convex Objectives}\label{sec:nonconvexEx}
\begin{figure}[htb!]
      \centering
      \includegraphics[scale=.25]{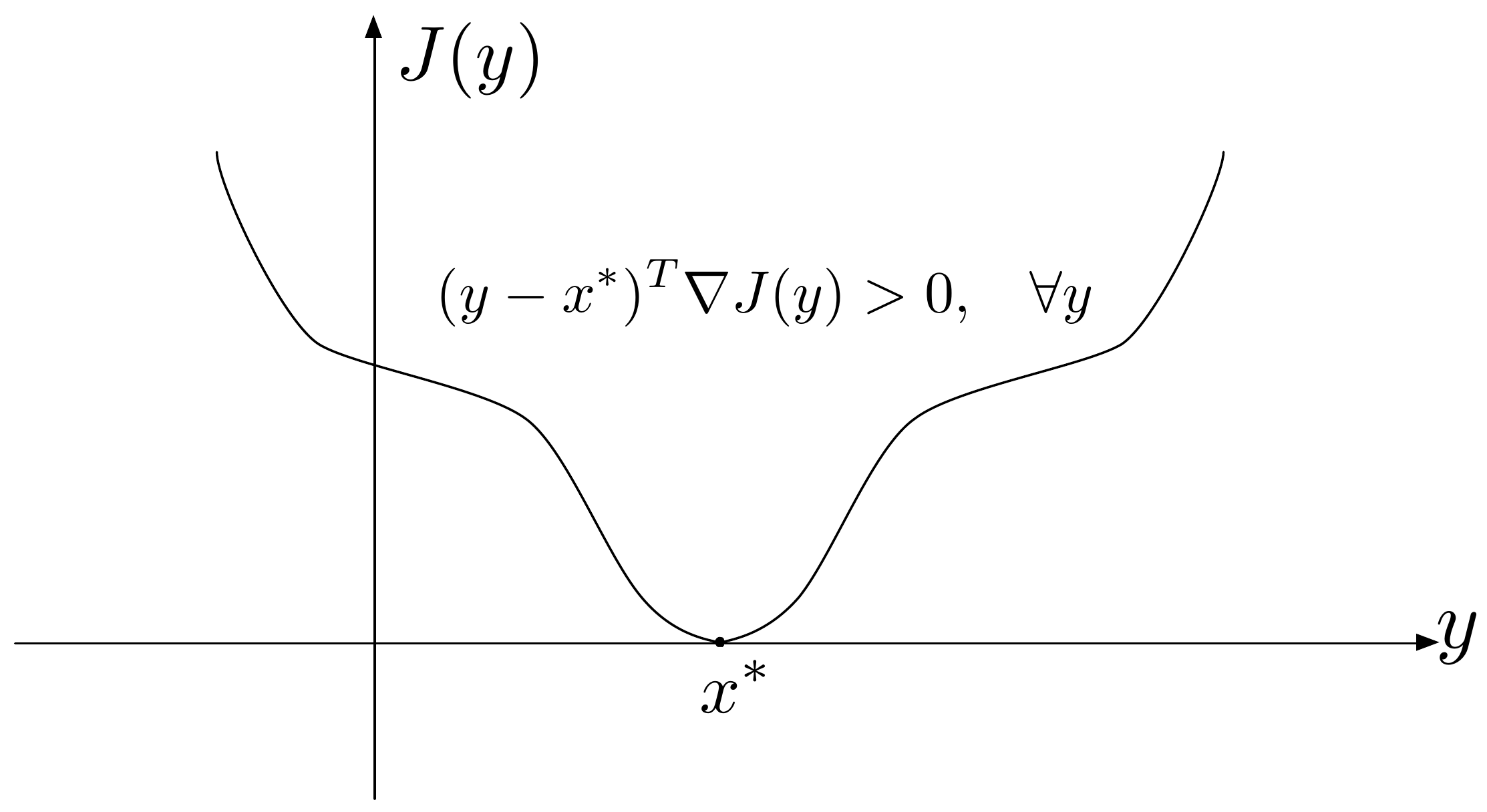}
      \caption{Gradients of non-convex, differentiable objective functions may satisfy the SP condition.}
      \label{fig:nonConvex}
      \vspace*{-.1in}
   \end{figure}

The examples in the prequel show how, in the context of optimization, the convexity of the objective function $J(\cdot)$ can be exploited to conclude that gradient-related seach directions $\Psi(\cdot)$ have the SPSP property. However, within the context of optimization, convexity is not necessary for SPSP, as shown in Figure \ref{fig:nonConvex}. The figure depicts a smooth, non-convex function $J(\cdot)$ for which $\Psi(\cdot)=\{\del J(\cdot)\}$ is SP with respect to $V(y)=\aHalf (y-x^*)^2$. 

This observation serves as another motivation for introducing the notion of SPSP search directions. Although we do not explore the possibility further here, the observation that convexity is not necessary for the SPSP property to hold may potentially be leveraged to inform novel designs of iterative optimization methods of the form \eqref{eq:iterative}. 

\section{SPSP Search Directions Generate SPAS Behavior}\label{sec:SPSP=SPAS}
A third motivation for introducing the SPSP property is that the qualitative behavior of sequences generated by numerical methods of the form \eqref{eq:iterative} with SPSP search directions is guaranteed to exhibit not only attractivity, but also stability with respect to the posited attractor $\mathcal{A}$. The concept of stable behavior is especially pertinent to applications involving the dynamic interaction between algorithms such as \eqref{eq:iterative} and physical systems. In such applications, the sequences $(y(t))_{t\in\mathbb{N}}$ may represent physical signals whose large excursions cannot be safely tolerated. 

In this section we investigate a number of scenarios under which the generic iterative methods of the form \eqref{eq:iterative}, having SPSP search directions $\Psi(\cdot)$, have an attractor $\mathcal{A}$ which is \emph{semiglobally, practically, asymptotically stable} (SPAS). We provide a precise definition of SPAS in Appendix \ref{sec:SPASdefs}, and Theorem \ref{thm:spasNew}, characterizes SPAS for a class of systems more general than those represented by \eqref{eq:iterative}. Our goal here is to link the SPSP property to the conditions of Theorem \ref{thm:spasNew}.

In the following two lemmas, we derive basic expressions for upper bounds on the difference $V(y^+)-V(y)$, under two separate sets of assumptions: either $V(\cdot)$ is the square of the distance to $\mathcal{A}$ and the constraint set $\Xi$ is generic, or, the form of $V(\cdot)$ is generic but \eqref{eq:iterative} does not involve a projection operation (i.e., $\Xi=\mathbb{R}^n$)\footnote{It remains an open problem to generalize the class of feasible sets for which similar inequalities can be derived for generic Lyapunov functions $V(\cdot)$.}.

From these upper bounds, we show in $\S$\ref{sec:conditions} how the conditions of Theorem \ref{thm:spasNew} can be satisfied when the search directions in \eqref{eq:iterative} possess the SPSP property and satisfy any one of three additional ``regularity'' assumptions, each of which is standard in the context of optimization.

\begin{lemma}[The case in which $V(y)=\aHalf \|y-\Pa(y)\|^2$ and $\Xi$ is generic]\label{lem:genericXi}
Consider the numerical method \eqref{eq:iterative}, and assume that $\Xi$ is a closed, convex set, though not necessarily bounded. Then, for all $y\in\Xi$,
\begin{equation}\label{eq:DeltaVGeneric1}
\Delta V(y) \leq -\alpha \del V(y)^Ts + \aHalf\alpha^2\|s\|^2,
\end{equation}
where $s\in\Psi(y)$, $V(y)=\aHalf\|y-\Pa (y)\|^2$, and $\mathcal{A}$ is some compact, convex subset of $\mathbb{R}^n$.  
\proof
We have that 
\begin{align}
\Delta V(y) &= \aHalf \|y^+-\Pa(y^+)\|^2 - \aHalf \|y-\Pa(y)\|^2
\nonumber \\
&\leq  \aHalf \|y^+-\Pa(y)\|^2 - \aHalf \|y-\Pa(y)\|^2
\label{eq:distIneq}\\
&= \aHalf \|\PXi(y-\alpha s(y)) - \Pa(y)\|^2 - \aHalf \|y-\Pa(y)\|^2
\nonumber \\
&= \aHalf \|\PXi(y-\alpha s(y)) - \PXi(\Pa(y))\|^2 - \aHalf \|y-\Pa(y)\|^2
\label{eq:XstSubsetXi} \\
&\leq \aHalf \|y-\alpha s(y)- \Pa(y)\|^2 - \aHalf \|y-\Pa(y)\|^2
\label{eq:PrNonExp}\\
&=- \alpha (y-\Pa(y))^Ts(y) + \aHalf\alpha^2 \|s(y)\|^2
\nonumber \\
&= -\alpha \del V(y)^Ts(y) + \aHalf\alpha^2 \|s(y)\|^2,
\label{eq:lastTime}
\end{align}
where \eqref{eq:PrNonExp} holds because the orthogonal projection operator is non-expansive, \eqref{eq:XstSubsetXi} holds because $\Xi\supset\mathcal{A}$ by assumption, and \eqref{eq:distIneq} holds because $\|y^+-\Pa(y^+)\|$ is the smallest distance between $y^+$ and $\mathcal{A}$, and therefore $\|y^+-\Pa(y^+)\|\leq \|y^+-\Pa(y)\|$, for any $y\in\mathbb{R}^n$.
$\diamondsuit$\end{lemma}

\begin{lemma}[The case in which $V(\cdot)$ is generic, and $\Xi=\mathbb{R}^n$]\label{lem:genericV}
Consider the numerical method \eqref{eq:iterative} with $\Xi=\mathbb{R}^n$ and a continuously differentiable function $V:\mathbb{R}^n\rightarrow\mathbb{R}$ having a locally Lipschitz gradient. For any compact $\Omega\subset\mathbb{R}^n$, there exists a number $L_{\del V}\in\mathbb{R}_{++}$ such that 
\begin{equation}\label{eq:DeltaVGeneric2}
\Delta V(y) \leq -\alpha\del V(y)^Ts + \alpha^2 L_{\del V}\|s\|^2,
\end{equation}
for all $y\in\Omega$ and $s\in\Psi(y)$. 
\proof
Using the Mean Value Theorem, it can be shown that 
\begin{align}\label{eq:initBd}
V(y^+) &\leq V(y) + \del V(y)^T(y^+-y) + L_{\del V}\|y^+-y\|^2,
\end{align}
where $L_{\del V}$ is the Lipschitz constant associated to $\del V(\cdot)$ on a set $\hat{\Omega}\supset\Omega$, constructed so that $y-\alpha s\in\hat{\Omega}$ whenever $y\in\Omega$ and $s\in\Psi(y)$. 

Since $y^+-y = -\alpha s$, \eqref{eq:initBd} yields the desired inequality \eqref{eq:DeltaVGeneric1}. 
$\diamondsuit$\end{lemma}

Since the inequalities \eqref{eq:DeltaVGeneric1} and \eqref{eq:DeltaVGeneric2} differ only by a constant factor in the last term, in the sequel we refer to the bound
\begin{align}\label{eq:DeltaVGeneric}
\Delta V(y) &\leq -\alpha \del V(y)^T s + \alpha^2 w\|s\|^2,
\end{align}
for some $w\in\mathbb{R}_{++}$.

\subsection{Meeting the Conditions of Theorem \ref{thm:spasNew}}\label{sec:conditions}

We now show how the conditions of Theorem \ref{thm:spasNew} can be satisfied when the search directions $s$ possess the SPSP property, in addition to satisfying any one of three conditions that are typically imposed in the optimization literature.

\begin{theorem}[$s$ satisfies a relative growth condition]
Consider the numerical method \eqref{eq:iterative}, and assume that $\Psi(\cdot)$ is SPSP (q.v. Definition \ref{def:SPStrictPseudo}) with respect to a function $V(\cdot)$, where $V(\cdot)$ and \eqref{eq:iterative} satisfy the conditions of either Lemma \ref{lem:genericXi} or \ref{lem:genericV}. Suppose that there exists a number $\beta\in\mathbb{R}_{++}$ such that for all $y\in\Xi$ and for all $s\in\Psi(y)$, $s$ satisfies the relative growth condition
\begin{equation}\label{eq:growthCon}
\|s\|^2\leq \beta\del V(y)^Ts.
\end{equation}
Then, $\mathcal{A}$ is SPAS (q.v. Definition \ref{def:spasNew2}) for \eqref{eq:iterative}, with Lyapunov function $V(\cdot)$. 
\proof
From \eqref{eq:DeltaVGeneric} and \eqref{eq:growthCon}, we have
\begin{equation*}
\Delta V(y) \leq -\alpha\del V(y)^Ts \big(1-\alpha w \beta\big).
\end{equation*}
Since $\Psi(\cdot)$ is SPSP with respect to $V(\cdot)$, for any desired $\sigma_o\in\mathbb{R}_{++}$, the following inequality holds:
\begin{equation*}
\Delta V(y) \leq
\begin{cases}
 -\alpha \phi_{\sigma_o,\epsilon}(y) \big(1-\alpha w \beta\big),
 & y\in \bar{B}_{\sigma_o}(\mathcal{A})\backslash B_{\epsilon}(\mathcal{A})
 \\
 \alpha b(1-\alpha w\beta),
 & y\in \bar{B}_{\epsilon}(\mathcal{A}),
\end{cases}
\end{equation*}
where $\phi_{\sigma_o,\epsilon}(\cdot)$ is continuous, positive on $\Xi\cap \big(\bar{B}_{\sigma_o}(\mathcal{A})\backslash B_{\epsilon}(\mathcal{A})\big)$ and radially unbounded with respect to $\mathcal{A}$ on $\Xi$. From this inequality we observe that the conditions of Theorem \ref{thm:spasNew} are satisfied with $\epsilon_o=\epsilon$, $\rho_o\equiv 0$,  $W_{\sigma_o,\epsilon_o}(y) = \alpha \phi_{\sigma_o,\epsilon_o}(y) \big(1-\alpha w \beta\big)$, $\pi=\alpha$, and $P_o = (0, \min\{ \tfrac{1}{w\beta},\alpha_1\}]$, where $\alpha_1$ is selected so that for any desired $b_o\in\mathbb{R}_{++}$, $\alpha_1 b(1-\alpha_1 w\beta)\leq b_o$.
$\square$\end{theorem}

\begin{rem}
A condition similar to the growth condition \eqref{eq:growthCon} is considered in $\S$2.2.3 of \cite{Polyak87}. 
$\diamondsuit$\end{rem}

\begin{theorem}[$s$ is locally bounded]
Consider the numerical method \eqref{eq:iterative}, and assume that $\Psi(\cdot)$ is SPSP with respect to a function $V(\cdot)$, where $V(\cdot)$ and \eqref{eq:iterative} satisfy the conditions of either Lemma \ref{lem:genericXi} or \ref{lem:genericV}. 

Suppose that for every $\sigma\in\mathbb{R}_{++}$, there exists a number $B\in\mathbb{R}_{++}$ such that for all $y\in\Xi\cap\bar{B}_{\sigma}(\mathcal{A})$ and for all $s\in\Psi(y)$, 
\begin{equation}\label{eq:sbddCon}
\|s\|\leq B.
\end{equation}
Then, $\mathcal{A}$ (q.v. Definition \ref{def:SPStrictPseudo}) is SPAS (q.v. Definition \ref{def:spasNew2}) for \eqref{eq:iterative}, with Lyapunov function $V(\cdot)$. 
\proof
Let $\sigma_o\in\mathbb{R}_{++}$ be arbitrary as in the statement of Theorem \ref{thm:spasNew}, and apply \eqref{eq:sbddCon} to \eqref{eq:DeltaVGeneric}, taking $\sigma=\sigma_o$. Then,
\begin{align}
\Delta V(y) &\leq -\alpha \del V(y)^Ts(y) + w\alpha^2B^2, \quad\forall y \in\bar{B}_{\sigma_o}(\mathcal{A}).
\end{align}
Since $\Psi(\cdot)$ is SPSP with respect to $V(\cdot)$ on $\Xi$, there exists a continuous function $\phi_{\sigma_o,\epsilon}(\cdot)$, which is positive on $\Xi\cap \big(\bar{B}_{\sigma_o}(\mathcal{A})\backslash B_{\epsilon}(\mathcal{A})\big)$ and radially unbounded with respect to $\mathcal{A}$ on $\Xi$, and a number $b\in\mathbb{R}_{++}$ such that
\begin{align}
\Delta V(y) &\leq
\begin{cases}
-\alpha \phi_{\sigma_o,\epsilon}(y) + w\alpha^2B^2
&y\in \bar{B}_{\sigma_o}(\mathcal{A})\backslash B_{\epsilon}(\mathcal{A}),
\\
\alpha b + w\alpha^2B^2
&y\in \bar{B}_{\epsilon}(\mathcal{A}).
\end{cases}
\end{align}

To show that $\Delta V(\cdot)$ satisfies the conditions of Theorem \ref{thm:spasNew}, take $\epsilon_o=\epsilon$, and any $b_o\in\mathbb{R}_{++}$ and  $\rho_o\in\mathbb{R}_{++}$ (but such that $\rho_o+\epsilon_o<\sigma_o$). Then, apply Lemma \ref{lem:underEstimation} to $\phi_{\sigma_o,\epsilon}(\cdot)$, with $\hat{\sigma}=\sigma_o$, $\hat{\epsilon}=\epsilon_o+\rho_o$, and taking $K_{\phi}$ to be any positive, real number satisfying
\begin{equation}
K_{\phi}> \frac{wB^2}{(\rho_o+\epsilon_o)^2}
\end{equation}
We thus obtain that
\begin{align}\label{eq:dVBdd}
\Delta V(y) &\leq
\begin{cases}
-\alpha^2 K_{\phi} \big(  \|y-\Pa(y)\|^2 - \tfrac{wB^2}{K_{\phi}}\big)
&y\in \bar{B}_{\sigma_o}(\mathcal{A})\backslash B_{\epsilon_o+\rho_o}(\mathcal{A}),
\\
\alpha b + w\alpha^2B^2
&y\in \bar{B}_{\epsilon_o+\rho_o}(\mathcal{A}),
\end{cases}
\end{align}
provided that
\begin{equation}\label{eq:alphaBoundLips}
\alpha\leq \alpha_q := \frac{c(\epsilon_o+\rho_o)^2}{\sigma_o^2 w B^2}, 
\end{equation}
with $c$ chosen such that $\Phi_c$ is a sublevel set of $\phi_{\sigma,\epsilon}(\cdot)$ that is strictly contained inside $\bar{B}_{\epsilon_o+\rho_o}(\mathcal{A})$ (q.v. Lemma \ref{lem:sublevelSetContainment}). 


From \eqref{eq:dVBdd} we see that all the conditions of Theorem \ref{thm:spasNew} are satisfied with 
\begin{equation}\label{eq:W_Lips}
W_{\sigma_o,\epsilon_o}(y) =\alpha^2 K_{\phi} \big(  \|y-\Pa(y)\|^2 - \tfrac{wB^2}{K_{\phi}}\big)
\end{equation}
and $P_o= (0,\alpha_o]$, where
\begin{equation*}
\alpha_o=\min\{\alpha_q,\alpha_1\},
\end{equation*}
and $\alpha_1$ is selected such that 
\begin{equation*}
\alpha_1 b +w\alpha_1^2 B^2 \leq b_o.
\end{equation*}
$\square$\end{theorem}

\begin{theorem}[$s$ is locally Lipschitz continuous]
Consider the numerical method \eqref{eq:iterative}, and assume that $\Psi(\cdot)$ is SPSP with respect to a function $V(\cdot)$, where $V(\cdot)$ and \eqref{eq:iterative} satisfy the conditions of either Lemma \ref{lem:genericXi} or \ref{lem:genericV}. Suppose that $\forall y\in\Xi$, $\Psi(y)$ is the singleton $\{s(y)\}$, and that for every $\sigma\in\mathbb{R}_{++}$ there exists a number $L\in\mathbb{R}_{++}$ such that for all $y_1,y_2\in\Xi\cap\bar{B}_{\sigma}(\mathcal{A})$, $\|s(y_1)-s(y_2)\|\leq L\|y_1-y_2\|$.
Then, $\mathcal{A}$ (q.v. Definition \ref{def:SPStrictPseudo}) is SPAS (q.v. Definition \ref{def:spasNew2}) for \eqref{eq:iterative}, with Lyapunov function $V(\cdot)$. 
\proof
Let $\sigma_o\in\mathbb{R}_{++}$ be arbitrary as in the statement of Theorem \ref{thm:spasNew}. From \eqref{eq:DeltaVGeneric} and the assumption that $\Psi(\cdot)$ is SPSP with respect to $V(\cdot)$ on $\Xi$, there exists a continuous function $\phi_{\sigma_o,\epsilon}(\cdot)$, which is positive on $\Xi\cap \big(\bar{B}_{\sigma_o}(\mathcal{A})\backslash B_{\epsilon}(\mathcal{A})\big)$ and radially unbounded with respect to $\mathcal{A}$ on $\Xi$, such that 
\begin{equation*}
\Delta V(y) \leq
-\alpha \phi_{\sigma_o,\epsilon}(y) + w\alpha^2\|s(y)\|^2, 
\end{equation*}
for all $y\in \bar{B}_{\sigma_o}(\mathcal{A})\backslash B_{\epsilon}(\mathcal{A})$. Using Young's inequality, we have that on the same set,
\begin{align}\label{eq:deltaVLips}
\Delta V(y) &\leq -\alpha \phi_{\sigma_o,\epsilon}(y) +w \alpha^2\|s(y)-s(\Pa(y))+s(\Pa(y))\|^2
\nonumber \\
&\leq -\alpha\phi_{\sigma_o,\epsilon}(y) + 2w \alpha^2\|s(y)-s(\Pa(y))\|^2 \nonumber\\
&\qquad 
+2w\alpha^2 s^*
\nonumber\\
&\leq -\alpha\phi_{\sigma_o,\epsilon}(y) + 2w \alpha^2L^2\|y-\Pa (y)\|^2 +2w\alpha^2 s^*,
\end{align}
where 
\begin{equation}
s^* = \max_{y\in \mathcal{A}}\|s(y)\|^2,
\end{equation}
exists because $\mathcal{A}$ is compact and $s(\cdot)$ is continuous.

To show that $\Delta V(\cdot)$ satisfies the conditions of Theorem \ref{thm:spasNew}, take $\epsilon_o=\epsilon$, and any $b_o\in\mathbb{R}_{++}$ and  $\rho_o\in\mathbb{R}_{++}$ (but such that $\rho_o+\epsilon_o<\sigma_o$). Then, apply Lemma \ref{lem:underEstimation} to $\phi_{\sigma_o,\epsilon}(\cdot)$, with $\hat{\sigma}=\sigma_o$, $\hat{\epsilon}=\epsilon_o+\rho_o$, and $K_{\phi}=2wL^2+\kappa$, where $\kappa$ is any real number satisfying
\begin{equation}\label{eq:kappaIneq}
\kappa > \frac{2ws^*}{(\epsilon_o+\rho_o)^2}.
\end{equation}
We thus obtain that
\begin{align}\label{eq:dVLips}
\Delta V(y) &\leq
\begin{cases}
-\alpha^2\kappa \big(  \|y-\Pa(y)\|^2 - \tfrac{2ws^*}{\kappa}\big),
\\
\quad\quad\quad y\in \bar{B}_{\sigma_o}(\mathcal{A})\backslash B_{\epsilon_o+\rho_o}(\mathcal{A}),
\\
\alpha b +2w\alpha^2 L^2 (\epsilon_o+\rho_o)^2 + 2w\alpha^2 s^*,
\\
\quad\quad\quad y\in
\bar{B}_{\epsilon_o+\rho_o}(\mathcal{A}),
\end{cases}
\end{align}
provided that
\begin{equation}\label{eq:alphaBoundLips}
\alpha\leq \alpha_q := \frac{c}{\sigma_o^2\big(2wL^2+\tfrac{2ws^*}{(\epsilon_o+\rho_o)^2}\big)}, 
\end{equation}
where $c$ is chosen such that $\Phi_c$ is a sublevel set of $\phi_{\sigma,\epsilon}(\cdot)$ that is strictly contained inside $\bar{B}_{\epsilon+\rho}(\mathcal{A})$ (q.v. Lemma \ref{lem:sublevelSetContainment}). 


From \eqref{eq:dVLips} we see that all the conditions of Theorem \ref{thm:spasNew} are satisfied with 
\begin{equation}\label{eq:W_Lips}
W_{\sigma_o,\epsilon_o}(y) =\alpha^2\kappa \big(  \|y-\Pa(y)\|^2 - \tfrac{2ws^*}{\kappa}\big)
\end{equation}
and $P_o= (0,\alpha_o]$, where
\begin{equation*}
\alpha_o=\min\{\alpha_q,\alpha_1\},
\end{equation*}
and $\alpha_1$ is selected such that 
\begin{equation*}
\alpha_1 b +2w\alpha_1^2 L^2 (\epsilon+\rho)^2 + 2w\alpha_1^2 s^* \leq b_o.
\end{equation*}
$\square$\end{theorem}

\begin{rem}[Related results for systems like \eqref{eq:iterative}]
For systems of the form $y^+\in y+\alpha F(y)$, where $F(\cdot)$ is a multifunction satisfying certain technical conditions, Theorem 2 in \cite{TeelCDC00a} can be applied to conclude a qualitative behavior similar to that described by Definition \ref{def:spasNew2}. However, in order to apply Theorem 2 from \cite{TeelCDC00a}, it is necessary to demonstrate that $\mathcal{A}$ is asymptotically stable for the system $\dot{y}\in F(y)$. 
$\diamondsuit$\end{rem}

\section{SPSP Implies Robustness}\label{sec:robustness}
Search directions $s$ for \eqref{eq:iterative}, selected from $\Psi(\cdot)$ having the SPSP property are \emph{robust} in the sense that the SPSP property is retained under sufficiently small relative and absolute deterministic errors. We formalize this observation in the following theorem. 

\begin{theorem}
Consider a  function $V\in C^1[\mathbb{R}^n,\mathbb{R}_+]$, which is positive definite and radially unbounded with respect to a compact set $\mathcal{A}\subset\mathbb{R}^n$, and has a locally Lipschitz gradient that is identically zero on $\mathcal{A}$. Consider also the multifunctions $\Psi:\mathbb{R}^n\rightrightarrows \mathbb{R}^n$ and 
\begin{equation}\label{eq:PsiHat}
\hat{\Psi}(y) = \{s+\eta(y),\; s\in\Psi(y)\},
\end{equation}
where $\eta(\cdot)$ has the property that for any $\sigma\in\mathbb{R}_{++}$, there exist positive, real numbers $a$ and $r$ such that
\begin{equation}\label{eq:errModel2}
\|\eta(y)\|\leq a + r\|y-\Pa(y)\|,
\end{equation}
for all $y\in\bar{B}_{\sigma}(\mathcal{A})$.

If $\Psi(\cdot)$ is SPSP with respect to $V(\cdot)$ on some $\Xi\supset \mathcal{A}$, then $\hat{\Psi}(\cdot)$ is also SPSP with respect to $V(\cdot)$ on $\Xi$, provided that for every $\sigma$, $a$ and $r$ are sufficiently small. 
\proof
To show that $\hat{\Psi}(\cdot)$ is SPSP with respect to $V(\cdot)$, we must demonstrate that for any $\hat{\sigma}\in\mathbb{R}_{++}$, and for some non-negative real numbers $\hat{\epsilon}$ and $\hat{b}$, there exists a function $\hat{\phi}_{\hat{\sigma},\hat{\epsilon}}\in C^0[\mathbb{R}^n,\mathbb{R}]$, which is positive on $\Xi\cap \big(\bar{B}_{\hat{\sigma}}(\mathcal{A})\backslash B_{\hat{\epsilon}}(\mathcal{A})\big)$ and radially unbounded with respect to $\mathcal{A}$ on $\Xi$, such that 
\begin{equation}
\del V(y)^T\hat{s} \geq
\begin{cases}
-\hat{b}, & \forall y\in\Xi\cap\bar{B}_{\hat{\epsilon}}(\mathcal{A}),\;\forall \hat{s}\in\hat{\Psi}(y),
\\
\hat{\phi}_{\hat{\sigma},\hat{\epsilon}}(y), & \forall y\in\Xi\cap (\bar{B}_{\hat{\sigma}}(\mathcal{A})\backslash B_{\hat{\epsilon}}(\mathcal{A})),\;\forall \hat{s}\in\hat{\Psi}(y),
\end{cases}
\end{equation}
whenever $a$ and $r$ are sufficiently small. 

Let $\hat{\sigma}$ be an arbitrarily large, positive, real number. By hypothesis, for any $\hat{\sigma}\in\mathbb{R}_{++}$ and for some non-negative numbers $\epsilon$ and $b$, there exists a function $\phi_{\hat{\sigma},\epsilon}\in C^0[\mathbb{R}^n,\mathbb{R}]$, which is positive on $\Xi\cap \big(\bar{B}_{\hat{\sigma}}(\mathcal{A})\backslash B_{\epsilon}(\mathcal{A})\big)$ and radially unbounded with respect to $\mathcal{A}$ on $\Xi$, such that 
\begin{equation}\label{eq:hyp}
\del V(y)^Ts \geq
\begin{cases}
-b, & \forall y\in\Xi\cap\bar{B}_{\epsilon}(\mathcal{A}),\;\forall s\in\Psi(y),
\\
\phi_{\hat{\sigma},\epsilon}(y), & \forall y\in\Xi\cap (\bar{B}_{\hat{\sigma}}(\mathcal{A})\backslash B_{\epsilon}(\mathcal{A})),\;\forall s\in\Psi(y). 
\end{cases}
\end{equation}

We aim to construct $\hat{b},\hat{\epsilon}$ and $\hat{\phi}_{\hat{\sigma},\hat{\epsilon}}(\cdot)$ from $b,\epsilon$ and $\phi_{\hat{\sigma},\epsilon}(\cdot)$. From \eqref{eq:PsiHat}, we have that for any $\hat{s}\in\hat{\Psi}(y)$,
\begin{equation}
\del V(y)^T\hat{s} = \del V(y)^T s + \del V(y)^T\eta(y),
\end{equation}
for some $s\in\Psi(y)$. Since $\del V(\cdot)$ is locally Lipschitz and identically zero on $\mathcal{A}$, 
\begin{align}\label{eq:niceBound}
\del V(y)^T\hat{s} &= \del V(y)^T s + \big(\del V(y)-\del V(\Pa (y))\big)^T\eta(y)
\nonumber \\
&\geq \del V(y)^T s - L_{\del V} \|y-\Pa (y)\| \big( a + r\|y-\Pa(y)\| \big),
\end{align}
for some $L_{\del V}\in\mathbb{R}_{++}$, and for all $y\in\bar{B}_{\hat{\sigma}}(\mathcal{A})$. By \eqref{eq:hyp} then, 
\begin{align}
\del V(y)^T\hat{s} &\geq  \phi_{\hat{\sigma},\epsilon}(y) - (a L_{\del V}) \|y-\Pa (y)\| - (r L_{\del V}) \|y-\Pa (y)\|^2,
\end{align}
for all $y\in\Xi\cap (\bar{B}_{\hat{\sigma}}(\mathcal{A})\backslash B_{\epsilon}(\mathcal{A}))$.

We apply Corollary \ref{cor:underEstimation} (q.v. Remark \ref{rem:doublebound}) to conclude that for any $\hat{\epsilon}\in (\epsilon,\hat{\sigma})$ there exists a number $c\in\mathbb{R}_{++}$ such that 
\begin{align}
\phi_{\sigma,\epsilon}(y) &\geq \frac{c}{2\hat{\sigma}}\|y-\Pa(y)\| + \frac{c}{2\hat{\sigma}^2}\|y-\Pa(y)\|^2 
\end{align}
for all $y\in\Xi\cap (\bar{B}_{\hat{\sigma}}(\mathcal{A})\backslash B_{\hat{\epsilon}}(\mathcal{A}))$. Therefore, on this set we have that $\del V(y)^T\hat{s}\geq \hat{\phi}_{\hat{\sigma},\hat{\epsilon}}(y)$, where 
\begin{align}
\hat{\phi}_{\hat{\sigma},\hat{\epsilon}}(y) = \big(\tfrac{c}{2\hat{\sigma}} - a L_{\del V}\big) \|y-\Pa(y)\| + \big(\tfrac{c}{2\hat{\sigma}^2}-r L_{\del V}\big)\|y-\Pa(y)\|^2
\end{align}
is positive provided that 
\begin{align}
a &< \frac{c}{2\hat{\sigma} L_{\del V}},\quad \text{and}
\\
r &< \frac{c}{2\hat{\sigma}^2 L_{\del V}}.
\end{align}
In these bounds, $c$ may be taken such that $\Phi_c$ is the largest sublevel set of $\phi_{\hat{\sigma},\epsilon}(\cdot)$ contained inside $\bar{B}_{\hat{\epsilon}}(\mathcal{A})$.


It remains to show that there exists a $\hat{b}\in\mathbb{R}_+$ such that $\del V(y)^T\hat{s}\geq-\hat{b}$, for all $y\in\bar{B}_{\hat{\epsilon}}(\mathcal{A})$. From \eqref{eq:hyp} we see that $\del V(y)^Ts$ is positive on $\bar{B}_{\hat{\epsilon}}(\mathcal{A})\backslash B_{\epsilon}(\mathcal{A})$, while it is no smaller than $-b$ on $\bar{B}_{\epsilon}(\mathcal{A})$. From this observation and inequality \eqref{eq:niceBound}, it follows that on $\bar{B}_{\hat{\epsilon}}(\mathcal{A})$, $\del V(y)^T\hat{s}\geq-\hat{b}$, where $\hat{b}= b + L_{\del V}\hat{\epsilon} (a + r\hat{\epsilon})$.
$\square$\end{theorem}

\section{Conclusions}
We introduced the notion of semiglobally, practically, strictly pseudogradient (SPSP) search directions for iterative numerical methods, and showed that a variety of optimization algorithms, including those affected by relative and absolute deterministic errors, have this property. We showed that iterative methods with SPSP search directions have semiglobally, practically, asymptotically stable attractors, and that the SPSP property is robust with respect to absolute and relative perturbations. Finally, we provided a set of technical lemmas that may serve as analytic tools to help establish the SPSP property in contexts other than those considered here.

Because iterative methods with SPSP search directions constitute a large class of optimization algorithms, and because the SPSP property implies SPAS, we anticipate that this property will be useful in guiding the design and analysis of novel data-driven management and control strategies for a variety of cyber-physical systems.

\appendix
\section{Semiglobal, Practical, Asymptotic Stability}\label{sec:SPASdefs}

We consider a class of discrete-time dynamical systems of the form
\begin{equation}\label{eq:system}
\xi^+ = \PXi\big[f(\xi;\pi)\big],\quad \xi\in\mathbb{R}^n,
\end{equation}
where $\Xi\subset\mathbb{R}^n$ is a closed, convex set, and $\pi\in\mathbb{R}^p$ parametrizes the function ${f:\mathbb{R}^n\rightarrow\mathbb{R}^n}$. 

We say that:

\begin{defn}\label{def:stabilityNew2}
A set $\mathcal{A}\subset\Xi$ is \emph{practically stable} for \eqref{eq:system} if for some $\check{\rho}_s\in\mathbb{R}_{++}$, and for any $\rho_s>\check{\rho}_s$, there exists a positive, real number $\delta$ and a set $P_s\subset\mathbb{R}^p$,  such that whenever $\pi\in P_s$ and $\xi(0)\in \bar{B}_{\delta}(\mathcal{A})\cap\Xi$, $\xi(t)\in\bar{B}_{\rho_s}(\mathcal{A})\cap\Xi$, for all $t\in\mathbb{N}$.
\end{defn}

\begin{defn}\label{def:attractivityGeneric}
A compact set $S\subset\mathbb{R}^n$ is  \emph{uniformly attractive for \eqref{eq:system} on a compact} $\Omega\subset\mathbb{R}^n$, if for every $\varepsilon\in\mathbb{R}_{++}$ for which $\bar{B}_{\varepsilon}(S)\cap\Xi\subset\Omega\cap\Xi$, there exists a number $T\in\mathbb{N}$ such that $\xi(t)\in\bar{B}_{\varepsilon}(S)$, whenever $\xi(0)\in\Omega$ and $t\geq T$.
\end{defn}

\begin{rem}
If a compact set $S\subset\mathbb{R}^n$ is  \emph{uniformly attractive} for \eqref{eq:system} on every compact $\Omega\subset\mathbb{R}^n$, then $S$ satisfies the usual definition of attractivity for \eqref{eq:system}. 
\end{rem}

\begin{defn}\label{def:attractivityNew2}
A compact set $\mathcal{A}\subset\Xi$ is  \emph{semiglobally, practically attractive} for \eqref{eq:system} if for some $\check{\rho}_a\in\mathbb{R}_{++}$, and for any $\sigma,\rho_a,\in\mathbb{R}_{++}$, with $\sigma>\rho_a>\check{\rho}_a$, there exists a set $P_a\subset\mathbb{R}^p$ such that whenever $\pi\in P_a$, the set $\bar{B}_{\rho_a}(\mathcal{A})$ is uniformly attractive for \eqref{eq:system} on $\bar{B}_{\sigma}(\mathcal{A})$.
\end{defn}

\begin{defn}\label{def:spasNew2}
A set $\mathcal{A}\subset\Xi$ is \emph{semiglobally practically asymptotically stable} (SPAS)  for \eqref{eq:system} if it is practically stable and semiglobally, practically attractive for \eqref{eq:system}.
\end{defn}

The following theorem characterizes the SPAS behavior defined above in terms of Lyapunov functions with certain properties. 

\begin{theorem}\label{thm:spasNew}
Consider the system \eqref{eq:system}, and suppose there exists a function $V\in C^0[\mathbb{R}^n, \mathbb{R}_+]$ which is radially unbounded and positive definite with respect to a compact set $\mathcal{A}\subset\Xi$ on $\mathbb{R}^n$. Suppose that for some $\epsilon_o\in\mathbb{R}_+$ and for any positive, real $\sigma_o$, $\rho_o$ and $b_o$ (with $\sigma_o>\epsilon_o+\rho_o$) 
there exists a set $P_o\subset\mathbb{R}^p$ and a function $W_{\sigma_o,\epsilon_o}\in C^0[\mathbb{R}^n,\mathbb{R}]$ such that whenever $\pi\in P_o$:
\begin{itemize}
\item P1: $W_{\sigma_o,\epsilon_o}(\xi) > 0$ for all  $\xi\in \Xi\cap (\bar{B}_{\sigma_o}(\mathcal{A})\backslash B_{\epsilon_o+\rho_o}(\mathcal{A}))$,
\item P2: $\Delta V(\xi)\leq -W_{\sigma_o,\epsilon_o}(\xi)$, for all $\xi\in\Xi\cap (\bar{B}_{\sigma_o}(\mathcal{A})\backslash B_{\epsilon_o+\rho_o}(\mathcal{A}))$, and 
\item P3: $\Delta V(\xi) \leq b_o$, for all $\xi\in \Xi\cap\bar{B}_{\epsilon_o+\rho_o}(\mathcal{A})$.
\end{itemize}

Then, $\mathcal{A}$ is SPAS for \eqref{eq:system}, with Lyapunov function $V(\cdot)$.
$\square$\end{theorem}

The proof is given in \cite{KvaternikSPAS}.


\section{Technical Lemmata}
The lemmas presented in this section provide analytic tools with a variety of applications in working with the SPSP property and the conditions of the SPAS Theorem \ref{thm:spasNew}. 

Lemma \ref{lem:sublevelSetContainment} states that for a continuous function that is radially unbounded with respect to some compact set $\mathcal{A}$, and positive on some ``band'' about that $\mathcal{A}$, one can always find a sublevel set of this function that fits inside an arbitrarily small ball containing $\mathcal{A}$. 

Lemma \ref{lem:underEstimation} and its corollary state that for any continuous function which is positive on a ``band'' surrounding some compact set, it is always possible to find either a linear or quadratic underestimator for the function on that band.

\subsection{Containment of Sublevel Sets}

\begin{lemma}\label{lem:sublevelSetContainment}
Consider a function $\phi_{\sigma,\epsilon}\in C^0[\mathbb{R}^n,\mathbb{R}]$ which is radially unbounded with respect to a compact set $\mathcal{A}\subset\mathbb{R}^n$, and positive on $\bar{B}_{\sigma}(\mathcal{A})\backslash \bar{B}_{\epsilon}(\mathcal{A})$, for some $\sigma\in\mathbb{R}_{++}\cup\{\infty\}$ and $\epsilon\geq 0$. For any  $\rho\in\mathbb{R}_{++}$ such that $\bar{B}_{\rho}(\bar{B}_{\epsilon}(\mathcal{A}))\subseteq\bar{B}_{\sigma}(\mathcal{A})$, there exists a number $l\in\mathbb{R}_{++}$ such that the set $\Phi_l=\{y\in\mathbb{R}^n\;|\;\phi_{\sigma,\epsilon}(y)\leq l\}$ is strictly contained
inside $\bar{B}_{\rho}(\bar{B}_{\epsilon}(\mathcal{A}))$.
%
\proof Figure \ref{fig:sublevelSetContainment} illustrates the constructions used in this proof. 
   \begin{figure}[htb!]
      \centering
      \includegraphics[scale=.25]{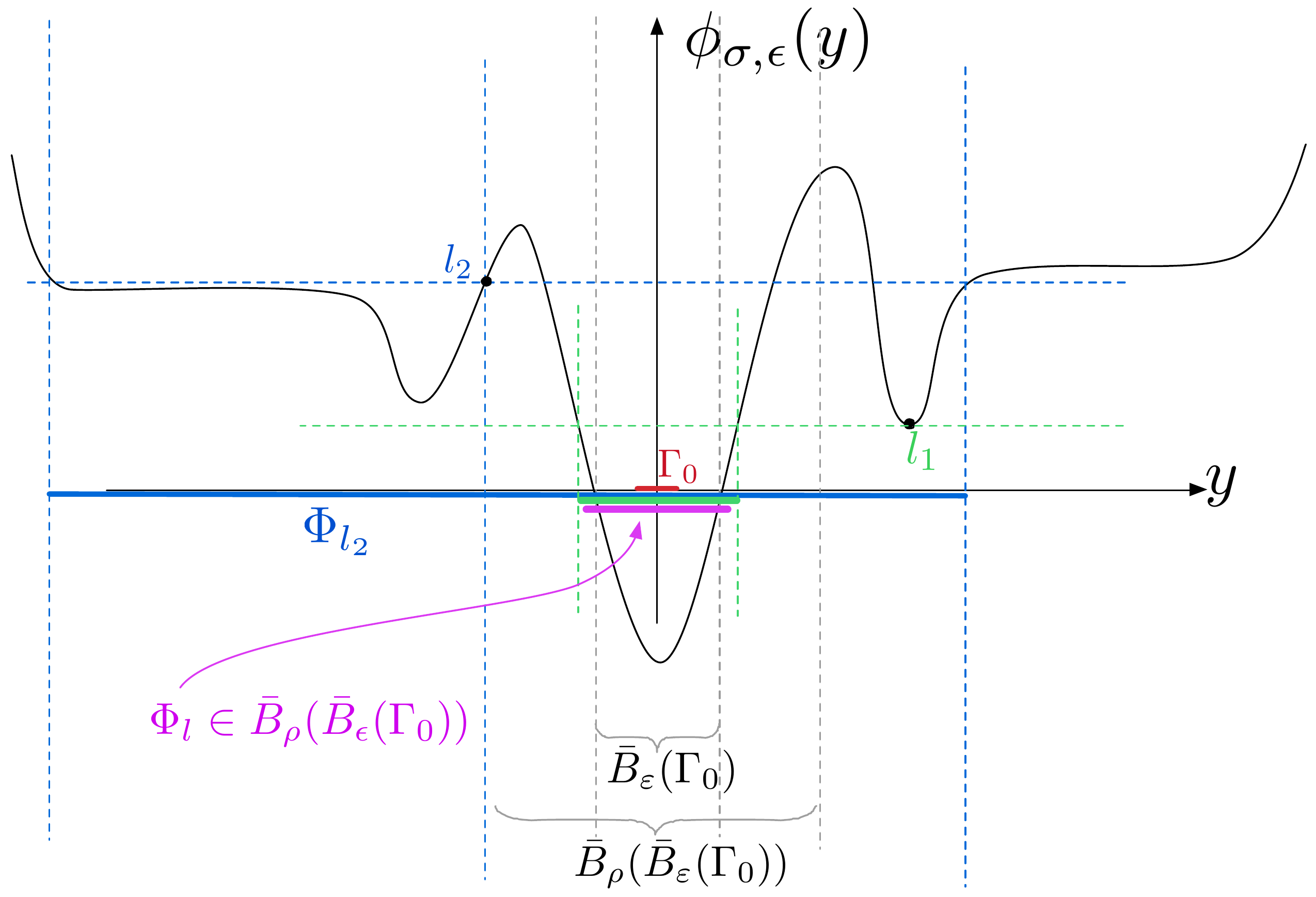}
      \caption{Constructions used in the proof of Lemma \ref{lem:sublevelSetContainment}.}
      \label{fig:sublevelSetContainment}
      \vspace*{-.1in}
   \end{figure}
Let $l_2$  be the minimum value attained by $\phi_{\sigma,\epsilon}(\cdot)$ on the set $\partial \bar{B}_{\rho}(\bar{B}_{\epsilon}(\mathcal{A}))$, and let $\Phi_{l_2}$ denote the $l_2$-sublevel set of $\phi_{\sigma,\epsilon}(\cdot)$. 
The existence of $l_2$ is guaranteed by the continuity of $\phi(\cdot)$ and the compactness of $\partial \bar{B}_{\rho}(\bar{B}_{\epsilon}(S))$. For the same reasons, the number 
$$
l_1 = \min_{y\in\Phi_{l_2}\backslash B_{\rho}(\bar{B}_{\epsilon}(\mathcal{A}))} \phi_{\sigma,\epsilon}(y)
$$
is also guaranteed to exist\footnote{The reason for performing the second minimization of $\phi_{\sigma,\epsilon}(\cdot)$ is illustrated in Figure \ref{fig:sublevelSetContainment}; the radial unboundedness and positivity of $\phi_{\sigma,\epsilon}(\cdot)$ do not preclude the possibility that some of the sublevel sets of $\phi_{\sigma,\epsilon}(\cdot)$ are not connected, since $\phi_{\sigma,\epsilon}(\cdot)$ is not required to be monotonically increasing in all directions away from $\mathcal{A}$. Therefore, it cannot be claimed that $\Phi_l\subset\bar{B}_{\rho}(\bar{B}_{\epsilon}(\mathcal{A}))$ for some $l\in(0,l_2)$. If all sublevel sets of $\phi_{\sigma,\epsilon}(\cdot)$ are connected, then $l_2\equiv l_1$, and the second minimization is superfluous.}.

Since $\phi_{\sigma,\epsilon}(\cdot)$ is positive on $\bar{B}_{\sigma}(\mathcal{A})\backslash\bar{B}_{\epsilon}(\mathcal{A})$, $l_1$ is positive and there exists a number $l\in(0,l_1)$. It can be seen that for any such $l$,  $\Phi_l$ is strictly contained inside $\bar{B}_{\rho}(\bar{B}_{\epsilon}(\mathcal{A}))$; otherwise, there would be a point $y_o\in \Phi_{l_2}\backslash B_{\rho}(\bar{B}_{\epsilon}(\mathcal{A}))$ for which both $\phi(y_o)\leq l$ and $\phi(y_o)\geq l_1$ hold. 
$\square$\end{lemma}

\begin{rem}
The statement of Lemma \ref{lem:sublevelSetContainment} can be contrasted with the immediate consequence of the definition of radial unboundedness, which provides the converse to this lemma: any sublevel set of a radially unbounded function is contained in a sufficiently large ball (and hence all its sublevel sets are bounded).
$\diamondsuit$\end{rem}

\begin{rem}
The behaivor of $\phi_{\sigma,\epsilon}(\cdot)$ inside $\bar{B}_{\epsilon}(\mathcal{A})$ is irrelevant to the conclusions of the Lemma; in particular, $\phi_{\sigma,\epsilon}(\cdot)$ need not be positive definite with respect to $\mathcal{A}$, and it need not be bounded below or above in $\bar{B}_{\epsilon}(\mathcal{A})$.
$\diamondsuit$\end{rem}

\begin{rem}
Lemma 2.5.1 in \cite{KvaternikPHD} can be regarded as a corollary of Lemma \ref{lem:sublevelSetContainment}. 
$\diamondsuit$\end{rem}

\subsection{Quadratic and Linear Underestimation on Bands}

\begin{lemma}\label{lem:underEstimation}
Consider a compact set $\mathcal{A}$, an arbitrary set $\Xi\supset \mathcal{A}$, and a function $\phi_{\sigma,\epsilon}\in C^0[\mathbb{R}^n,\mathbb{R}]$ which is positive on $\Xi\cap \big(\bar{B}_{\sigma}(\mathcal{A})\backslash B_{\epsilon}(\mathcal{A})\big)$, for some $\sigma>\epsilon\geq 0$. Then,
\begin{enumerate}
\item Given any three positive, real numbers $K_{\phi}$, $\hat{\sigma}$ and $\hat{\epsilon}$, with $\sigma\geq\hat{\sigma}>\hat{\epsilon}>\epsilon$, there exists a number $\alpha_q\in\mathbb{R}_{++}$ such that whenever $\alpha\in (0,\alpha_q]$,
\begin{equation}\label{eq:phiq}
\phi_{\sigma,\epsilon}(y)\geq\alpha K_{\phi}\|y-\mathbf{P}_{\mathcal{A}}(y)\|^2,\quad \forall y\in\bar{B}_{\hat{\sigma}}(\mathcal{A})\backslash B_{\hat{\epsilon}}(\mathcal{A}).
\end{equation}
\item Given any three positive, real numbers  $K_{\phi}$, $\hat{\sigma}$ and $\hat{\epsilon}$, with $\sigma\geq\hat{\sigma}>\hat{\epsilon}>\epsilon$, there exists a number $\alpha_l\in\mathbb{R}_{++}$ such that whenever $\alpha\in (0,\alpha_l]$,
\begin{equation}\label{eq:phil}
\phi_{\sigma,\epsilon}(y)\geq\alpha K_{\phi}\|y-\mathbf{P}_{\mathcal{A}}(y)\|,\quad \forall y\in\bar{B}_{\hat{\sigma}}(\mathcal{A})\backslash B_{\hat{\epsilon}}(\mathcal{A}).
\end{equation}
\end{enumerate}

\proof
Taking $\rho = \hat{\epsilon}-\epsilon$ and noting that $\bar{B}_{\rho+\epsilon}(\mathcal{A})\supseteq\bar{B}_{\rho}(\bar{B}_{\epsilon}(\mathcal{A}))$, we apply Lemma \ref{lem:sublevelSetContainment} to conclude that there exists a number $c\in\mathbb{R}_{++}$, such that the set
\begin{equation}\label{eq:Phicql}
\Phi_{c} = \{y\in\mathbb{R}^n\;\big|\; \phi_{\sigma,\epsilon}(y)\leq c\}
\end{equation}
is strictly contained inside $\bar{B}_{\hat{\epsilon}}(\mathcal{A})$.
%

Since $\|y-\mathbf{P}_{\!\mathcal{A}}(y)\|\leq\hat{\sigma}$ whenever $y\in\Bgo{\hat{\sigma}}$, and $\bar{B}_{\hat{\sigma}}(\mathcal{A})\supset \bar{B}_{\hat{\epsilon}}(\mathcal{A})\supset \Phi_c$, we have that 
\begin{equation}\label{eq:goodq}
\alpha K_{\phi}\|y-\mathbf{P}_{\!\mathcal{A}}(y)\|^2\leq \alpha_q K_{\phi}\hat{\sigma}^2,\quad \forall y\in\Bgo{\hat{\sigma}}\backslash\Phi_{c},
\end{equation}
for any positive, real $\alpha$ and $K_{\phi}$. On the same set, $\phi_{\sigma,\epsilon}(\cdot)$ is strictly larger than $c$. Therefore, taking
\begin{equation}\label{eq:alphaphiq}
\alpha\in(0,\alpha_q],\quad \alpha_q = \frac{c}{K_{\phi}\hat{\sigma}^2},
\end{equation}
implies that
\begin{equation}\label{eq:good2q}
\phi_{\sigma,\epsilon}(y)>c\geq \alpha_q K_{\phi}\hat{\sigma}^2\geq \alpha K_{\phi}\|y-\Pa(y)\|^2,
\end{equation}
whenever $y\in\Bgo{\hat{\sigma}}\backslash\Phi_{c}$. The bound \eqref{eq:phiq} follows since $\Bgo{\hat{\epsilon}}\supset\Phi_{c}$ by construction. 

Following the same arguments, the bound \eqref{eq:phil} can be seen to hold whenever $\alpha\in(0,\alpha_l]$, with 
\begin{equation}\label{eq:alphal}
\alpha_l = \frac{c}{K_{\phi}\hat{\sigma}}.
\end{equation}
$\square$\end{lemma}

\begin{corollary}\label{cor:underEstimation}
Consider a compact set $\mathcal{A}$, an arbitrary set $\Xi\supset \mathcal{A}$, and a function $\phi_{\sigma,\epsilon}\in C^0[\mathbb{R}^n,\mathbb{R}]$ which is positive on $\Xi\cap \big(\bar{B}_{\sigma}(\mathcal{A})\backslash \bar{B}_{\epsilon}(\mathcal{A})\big)$, for some $\sigma>\epsilon\geq 0$. Then, for any positive, real numbers $\hat{\sigma}$ and $\hat{\epsilon}$, with $\sigma\geq\hat{\sigma}>\hat{\epsilon}>\epsilon$, there exists a number $c\in\mathbb{R}_{++}$ such that 
\begin{equation}\label{eq:blank1}
\phi_{\sigma,\epsilon}(y)\geq \frac{c}{\hat{\sigma}^2}\|y-\mathbf{P}_{\mathcal{A}}(y)\|^2,\quad \forall y\in\bar{B}_{\hat{\sigma}}(\mathcal{A})\backslash B_{\hat{\epsilon}}(\mathcal{A}).
\end{equation}
and
\begin{equation}\label{eq:blank2}
\phi_{\sigma,\epsilon}(y)\geq \frac{c}{\hat{\sigma}}\|y-\mathbf{P}_{\mathcal{A}}(y)\|,\quad \forall y\in\bar{B}_{\hat{\sigma}}(\mathcal{A})\backslash B_{\hat{\epsilon}}(\mathcal{A}).
\end{equation}
\proof
Follows from the proof of Lemma \ref{lem:underEstimation} by constructing $c$ as in \eqref{eq:Phicql}, taking $K_{\phi}=1$, and setting $\alpha = \alpha_q$ to obtain \eqref{eq:blank1}, while setting $\alpha = \alpha_l$ to obtain \eqref{eq:blank2}, where $\alpha_q$ and $\alpha_l$ are as in \eqref{eq:alphaphiq} and \eqref{eq:alphal}, respectively.
$\square$\end{corollary}

\begin{rem}
If $\phi_{\sigma,\epsilon}(\cdot)$ is radially unbounded and positive definite with respect to $\mathcal{A}$, then we may take $\bar{B}_{\sigma}(\mathcal{A}) = \mathbb{R}^n$, and $\epsilon =0$.
$\diamondsuit$\end{rem}

\begin{rem}\label{rem:doublebound}
From \eqref{eq:blank1} and \eqref{eq:blank2} together, we can conclude that there exists a number $c\in\mathbb{R}_{++}$ such that 
\begin{equation}
\phi_{\sigma,\epsilon}(y) \geq \frac{c}{2\hat{\sigma}}\|y-\Pa(y)\| + \frac{c}{2\hat{\sigma}^2}\|y-\Pa(y)\|^2,\quad \forall y\in\bar{B}_{\hat{\sigma}}(\mathcal{A})\backslash B_{\hat{\epsilon}}(\mathcal{A}).
\end{equation}
$\diamondsuit$\end{rem}

\section*{Acknowledgments}
The author would like to thank Naomi Leonard for her keen insights and Biswadip Dey for some enjoyable technical discussions.

\newpage
 \bibliography{Bib}
\end{document}